\documentclass[a4paper, 11pt]{amsart}

\usepackage{amsmath}
\usepackage{amsfonts}
\usepackage{amssymb}
\usepackage{amsthm}
\usepackage{calc}
\def\geq{\geqslant}
\def\leq{\leqslant}
\def\Pic{\mathrm{pic}}
\def\Z{\mathbb{Z}}
\def\P{\mathbb{P}}
\def\deg{\mathrm{deg}\,}
\def\I{\mathcal{I}}
\def\Pic{\mathrm{Pic}} 
\def\Pl{(\P^n, \mbox{\rm line})}
\def\codim{\mathrm{codim}\hskip1pt}
\def\div{\mathrm{Div}}
\def\Supp{\mathop{\rm Supp}}
\def\invlim{\mathop{\mathrm{inv}\,{\lim}}}
\theoremstyle{definition}
\newtheorem{defin}{Definition}[section]

\newtheorem{example}[defin]{Example}
\theoremstyle{plain}              
\newtheorem{theorem}[defin]{Theorem}
\newtheorem{prop}[defin]{Proposition}
\newtheorem{lem}[defin]{Lemma}
\newtheorem{cor}[defin]{Corollary}

\newtheorem{question}[defin]{Question}
\newtheorem*{theoremA}{Theorem A}
\numberwithin{equation}{section}
\title[Birational geometry of rationally connected manifolds]{Birational geometry\\ of rationally connected manifolds\\ via quasi-lines}
\author{Paltin Ionescu}
\address{Department of Mathematics\\University of Bucharest\\ 14,
  Academiei Str.\\ RO--70109 Bucharest\\Romania}
\address{and}
\address{Institute of Mathematics of Romanian Academy\\P.O. Box
  1--764\\RO--70700 Bucharest\\Romania}
\email{Paltin.Ionescu@imar.ro}
\def\O{\mathcal{O}}
\def\model{(X,Y)}
\def\modelpr{(X',Y')}
\begin{document}

\begin{abstract}
This is, mostly, a survey of results about the birational geometry of rationally connected manifolds, using rational curves analogous to lines in ${\mathbb P}^n$ ({\it quasi-lines\/}). Various characterizations  of a Zariski neighbourhood of a line in ${\mathbb P}^n$ are obtained, some of them being new. Also, methods of formal geometry are applied for deducing results of birational nature. \end{abstract}

\subjclass[2000]{14-02, 14E08, 14E30, 14M20.}
\keywords{Rationally connected projective manifold, model, quasi-line, strongly rational, formal methods.}

\maketitle
\thispagestyle{empty}
\setcounter{section}{-1}
\section{Introduction}\label{s1}
The essential role played by rational curves in the birational
classification of algebraic varieties became clear since the
appearence of Mori theory \cite{14, 15}. In 1992, Koll\' ar, Miyaoka
and Mori \cite{12} introduced the very useful class of {\it rationally
  connected} manifolds. This class contains important subclasses, such
as unirational or Fano manifolds, admits several convenient
characterizations and has many good stability properties. One of the
characterizations of rationally connected manifolds is the presence of
a smooth rational curve having ample normal bundle.

A {\it model} of a rationally connected projective manifold $X$ is a
pair $(X,Y)$, where $Y$ is a smooth rational curve such that $N_{Y/X}$
is ample (cf.\ \cite{8}). Two models are equivalent, denoted
$(X,Y)\sim (X',Y')$, if there is an isomorphism between open Zariski
neighbourhoods of $Y$, respectively $Y'$, sending $Y$ to
$Y'$. Rationally connected manifolds are known to have a very
complicated birational geometry; see e.g.\ \cite{9, 10}. 
Hopefully, a convenient choice of models will simplify the original
birational problem. In fact, using Hironaka's result \cite{6}, we see
that the birational classification (of rationally connected manifolds)
is 
essentially the same as the classification of models, modulo the above
equivalence relation. Given a rationally connected projective manifold,
there is some birational model containing a ``quasi-line'', i.e.\ a
smooth rational  curve having the normal bundle of a line in $\P^n$
(cf.\ \cite{7}). It turns out that the Hilbert scheme of quasi-lines
has nice properties; in particular, counting curves through two
general points  can detect  the (uni)rationality of the given
manifold. An intermediate step in studying the equivalence of two
models $(X,Y)$ and $(X', Y')$ is the comparison of the formal
completions $X_{/Y}$ and $X'_{/Y'}$. This allows one to use powerful
results from formal geometry due to Hironaka, Hironaka--Matsumura,
Hartshorne, Gieseker (see \cite{1} or \cite {4}) in the study of
models.

This paper is, mostly, a survey based on the works \cite{2}, \cite{8} and
\cite{7}; on the other hand, it also contains several new results,
for instance Proposition~\ref{unu-opt}, Theorem~\ref{doi-opt}, Proposition~\ref{patru-doispe},
Theorem~\ref{patru-saptispe} and  Corollary~\ref{patru-optispe}. 
 Proposition~\ref{doi-unu}, Corollary~\ref{doi-trei}, Theorem~\ref{doi-cinci} and Proposition~\ref{cinci-*} are refinements of results from \cite{2} and \cite {7}. In the first section we discuss, informally,
models in general. In Theorem~\ref{unu-trei} we prove (cf.\ \cite {7}
and \cite{8}) a {\it reduction} result that replaces a given model by
one of lower dimension, in the presence of a suitable linear
system. We consider {\it deformations} of models and recall a
finiteness result from \cite{8}. Also, a notion of {\it minimality} is
briefly disscused. Next we recall from \cite{7} the existence of
quasi-lines and explain their use. The second section contains two
characterizations of the basic model $\Pl$. The first one,
Theorem~\ref{doi-cinci}, cf.\ \cite{7}, is in terms of curves only. The
second one, Theorem~\ref{doi-opt} which is new, is in terms of linear
systems of ``maximal dimension''. We note one consequence of
Theorem~\ref{doi-cinci}, cf.\ \cite{7}. Rational manifolds containing
``big'' open subsets of $\P^n$ (these are called {\it
  strongly-rational}, cf.~\cite{2}) are closed with respect to small
deformations. No analogous result is known for rational manifolds. In
the third section, following \cite{7}, we prove a characterization of
rationality in terms of quasi-lines and a useful result on the
ascent of rationality. We apply the last to give a simple proof
that del Pezzo manifolds of degree $\geq 4$ are rational. The last
section investigates the formal geometry of quasi-lines. We introduce,
following \cite{7}, a ``local number'' associated with a quasi-line; 
Theorem~\ref{patru-doi} relates global and local invariants of quasi-lines
and has many useful consequences. A computable property of quasi-lines,
called regularity, is introduced. A typical application is the
following (Corollary~\ref{patru-sapte}): if two models $(X,Y)$ and
$(X', Y')$, with $Y, Y'$  regular quasi-lines, are formally
equivalent (i.e.\ $X_{/Y}\simeq X' _{/Y'}$ as formal schemes), then
they are equivalent.
 In particular, we get (Proposition~\ref{cinci-*}): Let $X,X'\subset \P^{d+1}$, $d=3,4,5$, be smooth Fano threefolds of degree $d$ and let $Y\subset X$, $Y'\subset X'$ be general conics. If $X_{/Y}\simeq X'_{/Y'}$, then there exists an isomorphism $\varphi: X\to X'$ such that $\varphi (Y)=Y'$.
Theorem~\ref{patru-saptispe} contains a new characterization of the
formal completion $\P^n_{/{\rm line}}$. It may be used to get a new
description of the model $\Pl$, Corollary~\ref{patru-optispe}.
We work over the field $\mathbb{C}$  and we use the standard notation
and terminology of Algebraic Geometry, as in \cite{5}.
\section{Models of rationally connected manifolds}\label{S2}
Let $X$  denote a {\it complex projective manifold} of dimension $n$ (usually
$n\geq 2$). We say that $X$ is {\it rationally connected} (cf.\ \cite{12})  if
  there is a rational curve passing through any two given points of $X$. The
  following theorem summarizes some of the main properties of rationally
  connected manifolds.
\begin{theoremA} {\rm (cf.\ \cite{12, 3} or \cite{10})}
{\rm (i)} Unirational manifolds and Fano manifolds are rationally connected;

\hskip3pt{\rm (ii)} being rationally connected is a birational property and is
invariant under smooth deformation;

{\rm (iii)} rationally connected manifolds are simply connected and satisfy
\[ H^0 (X, \Omega_X ^{\otimes m})=0 \mbox{ for } m>0 \mbox{ and } H^i(X, {\O}_X)=0
  \mbox{ for } i>0;\]

\hskip.5pt{\rm (iv)} $X$ is rationally connected  if and only if there is a smooth
rational curve $Y\subset X$ such that the normal bundle of $Y$ in $X$ is
ample;

\hskip3.5pt{\rm (v)} if $n=2$,  $X$ is rationally connected if and only if $X$
is rational.
\end{theoremA}
\begin{defin} A {\it model} of a rationally connected manifold $X$ as above,
  is a pair $(X, Y)$, where $Y\subset X$ is a smooth rational curve with
  $N_{Y/X}$ ample.
\end{defin}

 Isomorphism of models is defined in the obvious way. More importantly, we
 introduce the following.
\begin{defin} Two models $(X,Y)$ and $(X', Y')$ are {\it equivalent}, denoted
 $(X, Y)\sim (X',Y')$, if there are Zariski open subsets $U\subseteq X$ and $U'
 \subseteq X'$, satisfying $Y\subset U$ and $Y' \subset U'$, and an
 isomorphism $\varphi: U \to U'$, such that $\varphi(Y)=Y'$.
\end{defin}

Note that in this case $X$ and $X'$ are birationally equivalent.
Conversely, assume that $X$ and $X'$ are birationally equivalent rationally
connected manifolds. Choose a birational map $\varphi:X \dasharrow X'$. By Hironaka's
result (cf.\ \cite{6}), we may find a composition of blowing-ups
$\sigma:\widetilde X \to X$ such that $\overline \varphi: = \varphi \circ
\sigma$ is a birational morphism. Now, we can choose a model $(X', Y')$ such
that $\overline\varphi{}^{-1}$ is an isomorphism along $Y'$; in other words,
by letting $\widetilde Y := \overline\varphi {}^{-1}(Y')$, the models
$(\widetilde X, \widetilde Y)$ and $(X', Y')$ are equivalent. Note, in
particular, that proving the rationality of a given rationally connected
manifold $X$, amounts to showing the equivalence $(\widetilde X, \widetilde
Y)\sim ( \mathbb{P}^n, \mbox{line})$, for a suitably chosen model of some blow-up
$\widetilde X$ of $X$.

To state our first result, we introduce some more notation. If $(X,Y)$ and
$(X',Y')$ are models and $\varphi :X\to X'$ is a morphism such that
$\varphi(Y)=Y'$, we write $\varphi :(X, Y)\to (X', Y')$ and call it a {\it
  morphism of models}. Given a model $(X, Y)$, we write $Y'\sim Y$ if $Y'$ is
a general deformation of $Y$; in particular, we can consider the new model
$(X, Y')$.
Recall that for any model $(X, Y)$, a theorem due to Grothendieck tells us that
$N_{Y/X}\simeq \bigoplus\limits_{i=1}^{n-1} \O_{\mathbb {P}^1} (a_i)$, $0<
a_1\leq \cdots\leq a_{n-1}$. Now we can state the following:
\begin{theorem}\label{unu-trei} {\rm (cf.\ \cite[Theorem 1.12]{8}, \cite[Theorem 1.7]{7})} Let $\model$ be a model such that $N_{Y/X}
  \simeq \bigoplus\limits_{i=1}^{n-1} \O_{Y}(a_i)$ with $a_1\leq \cdots \leq
  a_{n-1}$ and let $D$ be a divisor on $X$ such that  $(D \cdot Y):= d>0$,
  with $a_1\geq d$ and $\dim |D| \geq d$.    
Then there are $\widetilde X$ a blow-up of $X$, $Y'\sim Y$, and a diagram of
  models
\[(Z, \widetilde Y' ) \hookrightarrow (\widetilde X, \widetilde Y' )
  \overset{\varphi}{\longrightarrow} (\mathbb{P}^{\dim |D|-d+1}, l),\]
where $l$ is a line, such that: 

\hskip6pt{\rm (i)} $\varphi$ is surjective with connected fibres;

\hskip3pt{\rm (ii)} any smooth fibre of $\varphi$ is rationally connected;

{\rm (iii)} $Z:= \varphi^{-1} (l)$ is smooth;

\hskip.5pt{\rm (iv)} $\widetilde Y' $ is a section for $\varphi|_Z$.
\end{theorem}
\begin{proof}[Sketch of the proof] 
{\sc Step} 1. We may assume $|D|$ free of  fixed components and, replacing $Y$ by a
general deformation of it, we may also assume that $Y\cap
Bs|D|=\emptyset$. Next, we blow-up $d-1$ points on $Y$ and get a new manifold
$X'$. Replace the model $\model$ by the model $\modelpr$, where $Y'$ is the
proper transform of $Y$. We find the linear system $|D'|$ on $X'$ such that: 
$(D'\cdot Y')=1$ and $\dim |D'|\geq 1$.

{\sc Step} 2. By step 1, we may assume $d=1$ and $\dim |D|:= s\geq 1$. After suitable
blowing-up $\sigma :\widetilde X\to X$, we may also assume $Bs|D|=\emptyset$.
To prove (i), proceed by induction on $s$. If $s=1$, note that
 $(D \cdot \widetilde Y')=1$ shows that
the fibres of $\varphi_{|D|}$ are connected. If $s\geq 2$, $(D\cdot \widetilde
Y')=1$ implies that $\varphi$ is not composed with a pencil, so Bertini's
Theorem applies to find a smooth connected member $\Delta \in |D|$ through two
general points $x,y$ of $\widetilde X$. We may assume that $x,y\in \widetilde
Y' $ and we obtain the new model $(\Delta, \widetilde Y')$. The exact sequence
\[0\to \O_{\tilde X} \to \O _{\tilde X}(D) \to \O_{\Delta} (D)\to 0\]
shows that $\dim |D_{|\Delta}| = s-1$ and we may apply the induction
  hypothesis to deduce (i). 

(iii) follows by Bertini's Theorem, since the line determined by $\varphi (x),
  \varphi(y)$ is general.

(iv) is clear. 

To prove (ii), we use the model $(Z, \widetilde Y')$. We see that on a general
fibre of $\varphi_{|Z}: Z \to \mathbb{P}^1$,  two points may be joined by a
sequence of rational curves, making use of the same property on $Z$. 
  \end{proof}
\begin{cor}\label{unu-patru}{\rm (\cite[Corollary 1.8]{7})} Let $\model$ be a model and $D$ a divisor on $X$. Assume that
  $N_{Y/X} \simeq \bigoplus\limits_{i=1}^{n-1}\O _Y (a_i)$, $a_1\leq
  \cdots\leq a_{n-1}$, $d:= (D\cdot Y)>0$, $d\leq a_1$ and $\dim |D| \geq n+d
  -1$. Then $X$ is rational.
\end{cor}
\begin{cor} \label{unu-cinci} Keep all the other hypotheses of
  Corollary~{\rm \ref{unu-patru}} and
  replace the last inequality by $\dim |D|\geq n+d-2$. Then $X$ is birational
  to a conic bundle.
\end{cor}

Next, we consider deformations of models. Let $\model$ be a model. By a {\it
  deformation} of $\model$ we mean a commutative diagram
\[ \begin{array}{ccc}
\mathcal{Y}& \overset{i}{\longrightarrow}& \mathcal{X}\\[6pt]
\vcenter{\llap{$\scriptstyle q$}}
\Big\downarrow&\swarrow\!\vcenter{\rlap{$\scriptstyle p$}}&\\[10pt]
T\end{array}\]
where $p,q$ are proper smooth morphisms, $i$ is a closed embedding, $T$ is a
connected scheme such that $(\mathcal{X}_t, \mathcal{Y}_t)$ is a model for each (closed) $t\in T$
and $\model\simeq (\mathcal{X}_{t_0}, \mathcal{Y}_{t_0})$ for some $t_0\in T$. We say that
$\model$ and $(X', Y')$ are {\it deformation equivalent} if both appear as
fibres of the same deformation. By a {\it polarized model} we mean a triple
$(X,Y, H)$, where $\model$ is a model and $H$ is an ample divisor on $X$. Let
$d:= (H\cdot Y)$. Using Matsusaka's theorem (cf.\ \cite{13}) and its refinement
from \cite{11} together with Mori theory \cite{15}, one can prove the following
finiteness result:
\begin{theorem} {\rm (cf.\ \cite[Theorem 3.2]{8})} Fix $n\geq 2$ and $d> 0$. There are only
  finitely many isomorphism classes of polarized models $(X,Y, H)$ such that
  $\dim (X)=n$ and $(H\cdot Y)=d$, modulo deformations.
\end{theorem}

Next we discuss a notion of ``minimality'' for models, cf. \cite{8}.
\begin{defin}  A model $\model$ is {\it minimal} if any effective
    divisor $D$  satisfying $(D\cdot Y)=0$ is zero. 
\end{defin}

Note that for any model
    $\model$ the number of prime divisors $D$ with\linebreak $(D\cdot Y)= 0$ is finite.
In dimension two, given a model $\model$ there is a unique minimal model
    $(X_0, Y_0)$ and a birational morphism $\varphi: X \to X_0$ inducing an
    equivalence of $(X,Y)$ and $(X_0, Y_0)$ (here $X_0$ may be
    singular). Moreover, there is a complete classification of the pairs
    $(X_0, Y_0)$. See \cite[Proposition 1.21]{8} for details. In higher dimensions practically
    nothing is known about the existence (or uniqueness) of a minimal model in
    a given equivalence class of models.

The following proposition illustrates the use of the minimality property of a
given model.
\begin{prop}\label{unu-opt} Let $\varphi: X' \dasharrow X$ be a birational map inducing an
  equivalence of the models $(X',Y')$ and $(X,Y)$. Assume that $(X,Y)$ is
  minimal and for some (or any) ample line bundle $\mathcal{L} \in {\rm
  Pic}(X)$ there is $\mathcal{L}'\in \Pic (X')$ which is nef and coincides with $\varphi^*\mathcal{L}$ on the domain of $\varphi$. 
   Then $\varphi$ is a morphism.
\end{prop}

\begin{proof}
Let $\widetilde X$ be the normalization of the closure of the graph of $\varphi$, endowed with the natural projections $p:\widetilde X\to X'$ and $q: \widetilde X \to X$. We may assume $\mathcal{L}$ to be very ample. We have the equality of linear systems on $\widetilde X$,
$|p^*\mathcal {L}'|= E + |q^* \mathcal{L}|$ for some effective divisor $E$.
 We find easily that $(q_* E\cdot Y) = 0$; the minimality of $(X,Y)$ implies that $\codim_X(q_*E)\geq 2$. If $E>0$, it follows from Hodge Index Theorem, via suitable slicing, that there exists a curve $C\subset \widetilde X$ such that $q(C)$ is a point and $(C \cdot E)<0$. We infer that $(p^*\mathcal{L}'\cdot C)<0$ which contradicts the nefness of $\mathcal{L}'$. Thus we proved that $E=0$, so $p$ is an isomorphism.
\end{proof}
\begin{cor}\label{unu-noua} Let $\modelpr$ and $\model$ be equivalent minimal models. Assume,
  moreover, that $X$ and $X'$ are Fano manifolds. Then $\modelpr$ is isomorphic
  to $\model$.
\end{cor}

Next, we come to a very important question: Given $X$ a rationally connected
projective manifold, is there any ``convenient'' choice of a model $\model$?
As the model $(\mathbb{P}^n, \mbox{line})$ is the basic example, looking at
specific properties of a line in $\mathbb{P}^n$ suggests the following
definitions, cf.\ \cite{2}.
\begin{defin} (i) $Y\subset X$ is called a {\it quasi-line} if $N_{Y/X} \simeq
    \bigoplus\limits_{i=1}^{n-1}\O_{\mathbb{P}^1} (1)$;

(ii) $Y\subset X$ is called an {\it almost-line} if $Y$ is a quasi-line and
    there is a divisor $D$ on $X$ such that $(D\cdot Y)=1$.
\end{defin}
\begin{example}\label{unu-11}  Let $X$ be a Fano threefold of index two  such that
    $\Pic(X)=\Z[H]$, $H$ being the hyperplane section.
 If $Y \subset X$ is a general conic, $Y$ is a quasi-line. This was proved by
 Oxbury \cite{16}; see also \cite[Theorem 3.2]{2}
 for a more conceptual argument. It is
 easy to see that $X$ does not contain almost-lines.
\end{example}

The following theorem from \cite{7} is essential for the rest of this paper.
\begin{theorem}[existence of almost-lines] \label{unu-doispe} Let $X$ be a rationally connected projective manifold. There
  is $X'$ a smooth projective birational model of $X$, containing an
 almost-line.
\end{theorem}

Actually, one proves that $X'$ is got by a sequence of blowing-ups  of $X$ with
smooth two-codimensional centers; see \cite[Theorem 2.3, Proposition~2.5]{7}.

The importance of quasi-lines comes from the following considerations about the
Hilbert scheme of such curves. Let $Y\subset X$ be a quasi-line. The Hilbert
scheme of curves corresponding to the Hilbert polynomial (for a certain
polarization) of $Y$ in $X$ is smooth at $[Y]$. So $[Y]$ lies on a unique
irreducible component, say $\mathcal{H}$, of this Hilbert scheme. Hence, we can
speak about ``curves from the family determined by $Y$''. Note that this
applies to any smooth rational curve with ample normal bundle, and we
implicitly used it before, when speaking about $ Y' \sim Y$, a
general deformation of $Y$.
We have the universal family of such curves $\mathcal{Y}$ and the standard diagram
\[ \begin{array}{ccc}
\mathcal{Y}& \overset{\phi}{\longrightarrow}& X\\[6pt]
\vcenter{\llap{$\scriptstyle \pi$}}
\Big\downarrow&&\\[10pt]
\mathcal{H}&&\end{array}\]
Now fix a point $x\in Y$. Similar considerations apply to the closed subscheme
of $\mathcal{H}$ corresponding to curves of the family passing through $x$. We
denote by $\mathcal{H}_x$, respectively $\mathcal{Y}_x$, the Hilbert scheme of
these curves and their universal family. We keep the same notation $\pi$ and
$\phi$ for the restriction of the above projections to $\mathcal{Y}_x$.
\begin{defin} Let $Y\subset X$ be a quasi-line.

\hskip3pt
(i) The number of quasi-lines from the family determined by $Y$ passing
  through two general points of $X$ is denoted by $e(X,Y)$.

(ii) The number of quasi-lines from the family passing through one general
  point of $X$ and tangent to a general tangent vector at that point is
  denoted by $e_0(X,Y)$.
\end{defin}
It is easy to see that  given a model $\model$, $e(X,Y)$ and $e_0(X,Y)$ are
finite exactly when $Y$ is a quasi-line. Moreover, $e(X,Y)$ is nothing but the
degree of the projection $\phi:\mathcal{Y}_x \to X$ for a general point $x\in
X$. One can see also that we always have $e_0(X, Y)\leq e(X,Y)$. However, this
inequality may be strict.
\begin{example}\label{unu-paispe} (cf.\ \cite[Example 2.7]{2}) Take $X$ to be the desingularization  
 of the toric quotient of $\P^n $ $ (n\geq 3)$ by  a cyclic group of order $n+1$
 and $Y$ the quasi-line that is the image of a line in $\P^n$. We see that
 $e_0\model=e_0(\P^n, \mbox{line})=1$ and $e(X,Y)=n+1$.
\end{example}
\begin{example}\label{unu-cincispe} Let $X\subset \P^{d+1}$ be a smooth Fano threefold of degree $d$, $d=3,4,5$ and let $Y\subset X$ be a general conic (cf.\ Example~\ref{unu-11}). If 
$d=3$, we compute $e_0(X,Y)=e(X,Y)=6$, cf.\ \cite[Proposition 3.2]{7}. If $d=4$, we get similarly $e_0(X,Y)=e(X,Y)=2$. I thank Arnaud Beauville for kindly pointing out this fact to me. If $d=5$, we obtain $e_0(X,Y)=e(X,Y)=1$. 

\end{example}
\section{The model $(\P^n, {\mathop{\rm line}})$}\label{S3}
The main purpose of this section is to give two different characterizations of
models equivalent to the basic model $(\P^n, \mbox{line})$. A third
characterization will appear in Section~\ref{S5}, 
using formal geometry. The key to
all three characterizations is the following result, which is a
refinement of \cite[Theorem~4.4]{2}.
\begin{prop}\label{doi-unu} Let $\model$ be a model. The following
  conditions are equivalent:

\hskip3pt{\rm (i)} $\model \sim (\P^n, \mbox{\rm line})$;

{\rm (ii)} there is a divisor $D$ on $X$ such that $(D\cdot Y)=1$ and
 $\dim |D|\geq n$.
\end{prop}
\begin{proof} (i) $\Rightarrow$ (ii) is easy. 
To prove that (ii) $\Rightarrow $ (i), we first observe that $Y$ does not meet
the base locus of $|D|$. Otherwise, for some point $x\in Y$, we could find an
element in $|D|$, which is singular at $x$. But $Y$ deforms with the point
$x$ fixed, so replacing $Y$ by a general deformation of it passing through $x$,
we find that \mbox{$(D\cdot Y)\geq 2$}. The same argument shows two more
things. First, that $\dim |D|=n$; secondly, the rational map
$\varphi=\varphi_{|D|}$, which is defined along $Y$, is also \' etale along
$Y$. Moreover, $\varphi(Y)$ is a line in $\P^n$. It follows that the
restriction of $\varphi$ to $Y$ is an isomorphism onto a line and $Y$ is a
quasi-line. We may further assume, replacing $X$ by a suitable blowing-up,
that $\varphi:(X, Y)\to (\P^n, \mbox{line})$ is a morphism of models, \' etale
along $Y$. The following useful general lemma applies to show that $\varphi$
induces an equivalence $\model\sim (\P^n, \mbox{line})$, as
required. \end{proof}
\begin{lem} \label{doi-doi} Let $\varphi :(X', Y')\to (X, Y)$ be a morphism of models, with
  $Y, Y'$ quasi-lines, which is \' etale along $Y'$. Then $e(X, Y)=
  \deg\varphi \cdot e(X', Y')$.
\end{lem}
\begin{proof} Fix a point $y'\in Y'$ and let $y=\varphi(y')$. For $x\in X$ a
  general point, denote by $x'_1, \ldots, x'_d$ the points of the fibre
  $\varphi^{-1}(x)$, $d=\deg\varphi$. Consider the quasi-lines on $X'$
  equivalent to $Y'$, passing through $y'$ and some $x'_i$, $1\leq i\leq
  d$. Their images on $X$ are quasi-lines through $y$ and $x$, equivalent to
  $Y$. The induced map on Chow schemes $\varphi_* : {\rm Ch}_{y'} (X'
)\to
  {\rm Ch}_y (X)$ is injective when  restricted to the open sets parameterizing
  quasi-lines. This comes from the fact that the considered quasi-lines on
  $X'$ do not intersect the ramification divisor of $\varphi$. It follows
  that this restriction of $\varphi_*$ is also surjective, whence the desired
  equality. 
\end{proof}
 
\begin{cor} \label{doi-trei}
 Let $\model$ be a model such that there exists a divisor $D$ on $X$ with the
 following properties:

\hskip3pt{\rm (i)} $D$ is nef and big;

{\rm (ii)} $(D\cdot Y)=1$.

\noindent Then $\model \sim (\P^n, \mbox{\rm line})$.
\end{cor}
\begin{proof} 
Following \cite[p.\ 22]{2}: let $d:= (D^n)>0$. Consider the degree-$n$ Hilbert polynomial $p(t):= \chi (\O_X (tD))$. By duality and Kawamata--Viehweg vanishing theorem we have, for $i=1, \ldots, n$,
\[p(-i) =\chi (\O_X(-iD))= (-1)^n \chi (\O_X (K_X+iD))= (-1)^n h^0 (\O _X (K_X +iD)).\]
Since $(K_X \cdot Y) \leq -n-1$, we get $(K_X+ i D)\cdot Y<0$, hence $h^0 (\O_X (K_X +iD))=0$, $i=1, \ldots, n$. Therefore
\[ p(t)=\frac d {n!} (t+1)\cdots (t+n).\]
On the other hand,
\[p(-(n+2))=(-1)^n d(n+1) =(-1)^nh^0 (\O_X(K_X +(n+2)D)).\] 

Thus $h^0(\O_X(K_X+(n+2)D))=d (n+1) \geq n+1$.
 Put $D' := K_X+(n+2)D$. It follows that $\dim |D'|\geq n$ and $(D'\cdot Y)=1$, so Proposition~\ref{doi-unu} applies.
\end{proof}

\begin{cor}\label{doi-patru}
 Let $\model$ be a model such that there exists a divisor $D$ on $X$ with the
 following properties:

\hskip3pt{\rm (i)} $D$ is ample;

{\rm (ii)} $(D\cdot Y) =1 $.

\noindent Then $\model \overset{\sim}{\longrightarrow}  (\P^n, \mbox{\rm line}) $.
\end{cor}
 It is easy to see that  conditions (i) and (ii) imply that $(X,Y)$ is
 minimal, so the conclusion follows from Corollary~\ref{doi-trei} and Corollary~\ref{unu-noua}.
 
It is amusing to note that Corollary~\ref{doi-patru} (which follows also from
adjunction theory) implies the (well-known) fact that
$\P^n$ is closed with respect to (smooth) small deformations.

If $n=2$, $\model \sim \Pl$ means that there is a birational {\it morphism}
$\varphi :X\to \P^2$ such that $Y$ is the pull-back of a line. For any $n\geq
3$, there are models $\model \sim \Pl$ such that there is no birational
morphism $\varphi:X\to \P^n$.  In fact, there are models $(X,Y)\sim \Pl$ such that there is no divisor $D$ on $X$ which is nef and big with $(D\cdot Y)=1$; cf.\ \cite[Example~(4.7)]{2},  and \cite[Example (2.8)]{8}. 
  
  Now, we can state the first characterization of the models
  $\model$ equivalent to  $\Pl$. It uses only properties of the family of
  curves determined by $Y$. The result is a more precise form of Theorem~4.2 of \cite{7}.
\begin{theorem}\label{doi-cinci}  Let $\model$ be a
 model. The following assertions are equivalent:

{\rm (a)} $\model\sim \Pl$;

{\rm (b)} {\rm (i)} $Y$ is a quasi-line,

\hskip16pt{\rm (ii)} $e(X, Y)=1$, and

\hskip13pt{\rm (iii)} there exists a point $x\in Y$ such that, for any 
$Y'\sim Y$ with $x\in
Y'$, it follows that $x$ is a smooth point of $Y'$.
\end{theorem}
 
We shall see in the next section that (b)(i) and (b)(ii) together
imply that $X$ is rational. However, there are examples showing that condition
(b)(iii) is essential for the validity of  Theorem~\ref{doi-cinci} 
(see Example~\ref{patru-douazeci}).
\begin{proof}
Following \cite[Theorem 4.2]{7}, we sketch a proof that (b)(i), (b)(ii) and (b)(iii)
 together imply that $\model\sim \Pl$, the other implication being easier.
 We recall the basic diagram 
\[ \begin{array}{ccc}
\mathcal{Y}_x& \overset{\phi}{\longrightarrow}& X\\[6pt]
\vcenter{\llap{$\scriptstyle \pi$}}
\Big\downarrow&&\\[10pt]
\mathcal{H}_x&&\end{array}\]
where $\mathcal{H}_x$ is the Hilbert scheme of curves from the family determined
by $Y$, passing through $x\in Y$. $\pi$ has a section $\mathcal {E}$; condition
(b)(iii) implies that $\mathcal{E}$ is an effective Cartier divisor on
$\mathcal{Y}_x$ (see \cite[Lemma 4.3]{7}). Let $\sigma :X' \to X$ be the
blowing-up of $x\in X$, with $E\subset X'$ its exceptional divisor. We get a
birational morphism $\phi':\mathcal{Y}_x\to X'$ such that $\phi=\sigma
\circ\phi'$. Next we remark that, if $F$ is a general fibre of $\pi$ and $y\in
F\setminus (F\cap \mathcal{E})$, then $\phi$ is a local isomorphism at $y$. We
take $H \subset E\simeq \P^{n-1}$ a general hyperplane; we look at $D:= \phi
((\phi'_{|\mathcal{E}}\circ s \circ \pi)^*H)$, where $s: \mathcal {H}_x \to
\mathcal{E}$ is the inverse of the isomorphism $\pi_{|\mathcal {E}}: \mathcal{E}
\to \mathcal{H}_x$. We find that $x$ is a smooth point of the effective divisor $D$. A general deformation of $Y$ through $x$ meets $D$ only at $x$ and the intersection is transverse, so $(D \cdot Y)=1$. Next, for a fixed point $x\in Y$ satisfying (b)(iii), let us denote by $|D_x|$ the linear system constructed above.
We have $\dim |D_x|\geq n-1$. Either $\dim |D_x|\geq n$ or $\dim |D_x|=n-1$ and $x\in Bs|D_x|$. In the first case, we may apply Proposition~\ref{doi-unu} to conclude. If the last case happens for any $x\in Y$ satisfying the {\it open} condition (b)(iii), it follows that $\Pic(X)$ is uncountable. But $X$ is rationally connected, so $H^1(X, \O_X)=0$ and $\Pic (X)$ has to be countable. This is a contradiction. \end{proof}  
 
The following  definition first appeared in \cite{2}.
\begin{defin} A projective manifold $X$ is {\it strongly rational}, if there
    is a model $\model$ such that $\model\sim \Pl$. Equivalently, $X$ contains
    an open subset $U$ which is isomorphic to an open subset $V\subseteq \P^n$
    such that $\codim _{\P^n}(\P^n\setminus V)\geq 2$.
\end{defin} 
\begin{theorem} \label{doi-sapte}{\rm (cf.\ \cite[Proposition 3.10]{2}, \cite[Theorem 4.5]{7})} Let $X$ be a projective manifold. The
    following properties of $X$ are closed with respect to (smooth) small
    deformations:

\hskip6pt{\rm (i)} $X$ contains a quasi-line (respectively an almost-line);

\hskip3pt{\rm (ii)} $X$ contains a quasi-line $Y$ and $e(X,Y)=1$;

{\rm (iii)} $X$ is strongly rational.
\end{theorem}

The proof of (iii) follows from the characterization of strongly rational manifolds
given in Theorem~\ref{doi-cinci}.

Regarding the last point in Theorem~\ref{doi-sapte}, we recall that it is an
open problem whether or not small deformations of rational manifolds remain
rational. The example of cubic fourfolds in $\P^5$ seems to suggest that the
answer should be negative. The positive result in Theorem~\ref{doi-sapte}(iii)
illustrates how the use of models may simplify problems of birational nature.

The second characterization of models $\model$ which are equivalent to\linebreak $\Pl$
is via an extremality property of linear systems. It has been conjectured in
\cite[Conjecture 2.3]{8}.
\begin{theorem}\label{doi-opt} Let $\model$ be a model and $D$ a divisor on $X$.
Let $d=(D\cdot Y)$. We have:

\hskip3pt{\rm (i)} $\dim |D| \leq \binom {n+d}n -1$;

{\rm (ii)} $\model \sim \Pl$ if and only if equality holds in {\rm (i)}, for
some linear system $|D|$ with $d>0$
\end{theorem}
\begin{proof}
(i) (cf.\ \cite [Proposition~2.1]{8}) We may assume $d>0$. Consider the d-th jet
bundle of $\O_X(D)$, denoted $\mathcal{J}_d(D)$. Consider also the natural map $u: H^0 (X, \O_X (D)) \otimes \O_X\to \mathcal {J}_d (D)$ which sends a section to its d-th jet. We claim that $u$ is generically injective, which implies the desired inequality.

Let $s\in H^0(X, \O _X (D))$ be a non-zero section and let $D':= (s)_0$. If $x\in Y \cap \Supp (D')$, we may deform $Y$ to $Y'$ by keeping $x$ fixed, so that the local intersection number $(Y'\cdot D')_x$ be defined. If $u_x(s)=0$, we get $(Y'\cdot D')_x>d=(Y\cdot D)$, so $s$ has to be zero. Thus $u$ is injective in the open set swept out by the deformations of $Y$ having ample normal bundle. 

(ii) Assume that $\dim |D| = \binom {n+d} n -1 $. Denote 
by $\mu_x(D)$ the multiplicity at a point $x\in X$ of the effective divisor $D$. For $x\in Y$, let $\Lambda_x \subseteq |D|$ be the linear system 
\[\Lambda _x = \{ D' \in |D|\mid \mu_x(D') \geq d\}= \{ D' \in |D| \mid \mu _x (D') = d\}.\]
We have 
\[ \dim \Lambda _x \geq \binom {n+d} d - \binom {n+d- 1} { d-1} -1= \binom {n+ d-1}d - 1.\]

Let $\sigma : \widetilde X \to X$ be the blowing-up of $X$ at $x$ and let $E$ be its exceptional divisor. Write $\widetilde D $ for the proper transform of $D\in \Lambda _x$. Consider the exact sequence
\[ 0 \to \O_{\widetilde X} ( \sigma ^* (D)- (d+1) E)\to \O_{\widetilde X} (\widetilde D) \to \O_{\P^{n-1}} (d) \to 0.\]
We have $H^0(\widetilde X, \O_{\widetilde X}( \sigma ^* (D) - (d+ 1) E ))=0$, so $\dim |\widetilde D|\leq \dim |\O_{\P^{n-1}} (d)|= \binom {n+d -1} d - 1 $; but $\dim |\widetilde D| = \dim \Lambda _x \geq \binom{n+d-1} d -1$. So $|\widetilde D|$ has no base-points along $E$. If $\widetilde Y$ is the proper transform of $Y$, we also find: $(\widetilde D\cdot \widetilde Y)=0$. It follows that we have a rational map $\varphi := \varphi _{|\widetilde D|} : \widetilde X \dasharrow \P^N$, whose image is of dimension $n-1$ (because its restriction to $E$ is finite onto $\varphi (E)$) and contracts the proper transform $\widetilde Y'$ of a general deformation of $Y$. By generic smoothness, there is only one such curve $\widetilde Y'$ passing through the general point of $\widetilde X$, which means exactly that $e(X,Y)=1$. Note, in particular, that $Y$ has  to be a quasi-line. So we have checked the first two conditions in Theorem~\ref{doi-cinci}(b). Let us verify the third, too.  We may assume, after suitable blowing-up, that $|\widetilde D|$ is base-points free, see  \cite[p.\ 1068]{7}. Let $Y'\sim Y$, $x \in Y'$ and consider its proper transform, $\widetilde Y'$. We may find an element $\widetilde D' \in |\widetilde D|$ such that no irreducible component of $\widetilde Y'$ is contained in $\Supp (\widetilde D')$. It follows that there is some $D'\in \Lambda_x$ such that $(D' \cdot Y')_x $ is defined. But we have $\mu_x (D')=d$ which forces $Y'$ to be smooth at $x$. The conclusion follows now by applying Theorem~\ref{doi-cinci}.\end{proof}  

 Note that, for $d=1$, Theorem \ref{doi-opt} reduces to Proposition~\ref{doi-unu}
and was used in the proof, via Theorem~\ref{doi-cinci}. 

The following corollary is a new characterization of $\P^n$ (or of Veronese
varieties) in terms of models.
\begin{cor} \label{doi-noua}{\rm (cf.\ \cite{8})} Let $\model$ be a minimal model. The following
 are equivalent:

\hskip3pt{\rm (i)} $\model \simeq \Pl$;

{\rm (ii)} there is a divisor $D$ on $X$ such that $d=(D\cdot Y)>0$ and $\dim
|D|\geq \binom {n+d}n -1$.
\end{cor}
A different proof of Corollary~\ref{doi-noua}, based on Mori's characterization
of $\P^n$ as the only projective manifold with ample tangent bundle \cite{14},
was given in \cite[Proposition 2.10]{8}.
\section{(Uni)rationality via quasi-lines}\label{S4}
 The next proposition, taken from \cite{7}, shows how one may use quasi-lines for detecting the (uni)rationality of a given rationally connected projective manifold.

\begin{prop}\label{trei-trei} {\rm (cf.\ \cite[Proposition 3.1]{7})}
 Let $(X, Y)$ be a model with $Y$ a quasi-line.

\hskip3pt{\rm (i)} If $e_0(X,Y)=1$, then $X$ is unirational.

{\rm (ii)} If $e(X,Y)=1$, then $X$ is rational.
\end{prop} 
\begin{proof}
Let $x\in Y$ be a fixed point. Consider the standard diagram
\[ \begin{array}{ccc}
\mathcal{Y}_x& \overset{\phi}{\longrightarrow}& X\\
\vcenter{\llap{$\scriptstyle \pi$}}
\Big\downarrow&&\\
\mathcal{H}_x&&\end{array}\]
given by the family of curves determined by $Y$ and passing through $x$. $\pi$
has a section $\mathcal{E}= \phi^{-1}(x)$. Consider also $\sigma :{\rm Bl}_x(X)
\to X$, the blow-up of $X$ at $x$ and let $E$ be its exceptional divisor.
The rational map $\sigma^{-1} \circ \phi$ is defined at a general point of
$\mathcal{E}$ and maps $\mathcal{E}$ to $E$. The condition $e_0(X,Y)=1$ means
exactly that the restriction of the rational map $\sigma^{-1}\circ \phi$ to
$\mathcal{E}$ gives a birational isomorphism to $E$. But $\mathcal {Y}_x$ is
birationally a conic bundle for which $\mathcal{E}$ provides a rational
section. So $\mathcal{Y}_x$ is birational to $\mathcal{E}\times \P^1$ and (i)
follows.

To see (ii), observe that  the restriction map $\sigma^{-1} \circ \phi:
\mathcal{E}\dasharrow E$ is dominant, being generically finite. Since $e(X,Y)=1$,
$\sigma^{-1}\circ \phi$ is birational and by Zariski Main Theorem, it follows
that its restriction to $\mathcal {E}$ is also birational. Thus
$\mathcal{Y}_x$ is rational and so is $X$.\end{proof}  
\begin{defin}\label{trei-patru} Let $X$ be a rationally connected projective manifold. We denote
by $e(X)$ the minimum of $e(X',Y')$ for all models $(X',Y')$, where
$\sigma:X'\to X$ is a composition of blowing-ups with smooth centers and $Y'$
is a quasi-line in $X'$.
\end{defin}  
\begin{theorem}\label{trei-cinci} {\rm (cf.\ \cite[Theorem 3.4]{7})} Let $X$ be a rationally connected projective manifold.

\hskip3pt{\rm (i)} $e(X)$ is a birational invariant of $X$;

{\rm (ii)} $X$ is rational if and only if $e(X)=1$.
\end{theorem}
\begin{proof} 
 {\rm (i)} Let $\varphi : X_1 \dasharrow X_2$ be a birational isomorphism between two rationally connected projective manifolds. Let $\sigma : X' \to X_2$ be a composition of blowing-ups and let $Y' \subset X'$ be a quasi-line such that $e(X_2) = e(X',Y')$.
Let $\mu = \sigma^{-1} \circ \varphi:X_1 \dasharrow X'$. Take $\rho : X \to X_1$, a composition of blowing-ups such that $\mu \circ \rho : X \to X'$ is a birational morphism.
Let $Y \subset X$ be the inverse image by $\mu \circ \rho$ of a general deformation of $Y'$. It follows $e(X_2) = e(X',Y') = e(X,Y) \geq e(X_1)$. The opposite inequality follows by symmetry.

  {\rm (ii)} $e(\P^n) = 1$, so $e(X) = 1$ if $X$ is rational, by (i). The converse follows from Proposition~\ref{trei-trei} .
\end{proof}

The following result on the ascent of rationality was found in the context of 
Theorem~\ref{unu-trei}.
\begin{theorem}\label{trei-sase} {\rm (cf.\ \cite[Theorem 1.3]{7})}
Let $X$ be a projective variety and $|D|$ a complete linear system of Cartier divisors on it. Let $D_1,\ldots,D_s \in |D|$ and put $W_i := D_1\cap \cdots \cap D_i$ for $1\leq i\leq s$. Assume that $W_i$ is smooth, irreducible of dimension $n-i$, for all $i$. Assume moreover that there is a divisor $E$ on $W := W_s$ and a linear system $\Lambda \subset |E|$ such that:

\hskip3pt{\rm (i)} $\varphi_{\Lambda} : W \dasharrow \P^{n-s}$ is birational, and

{\rm (ii)} $|D_{|W} - E|\neq \emptyset $.
\noindent Then $X$ is rational.
\end{theorem}
\begin{proof}
  Induction on $s$. We explain the case $s=1$, the general case being
  completely similar. So, let $W \in |D|$ be a smooth, irreducible Cartier
  divisor such that $\varphi_{\Lambda} : W \dasharrow \P^{n-1}$ is birational for
  $\Lambda \subset |E|$, $E \in \div(W)$ and $|D_{|W} - E|\neq
  \emptyset$. Replacing $X$ by its desingularization, we may assume that $X$
  is smooth. As $W$ is rational, it is rationally connected, so we may find
  some smooth rational curve $Y \subset W$ with $N_{Y/W}$ ample. We have $(Y\cdot E) > 0$ and from (ii) we deduce $(Y \cdot D) > 0$. From the exact sequence of normal bundles we get that $N_{Y/X}$ is ample, so $X$ is rationally connected. In particular, $H^1(X, \O _{X}) = 0$. 
  
  The exact sequence 
  $$0 \to \O_X \to \O_X (D) \to \O_W (D) \to 0,$$ 
  shows that $\dim |D| = \dim |D_{|W}| + 1 \geq \dim |E| + 1 \geq n$.
  
  We may choose a pencil $(W, W') \subset |D|$, containing $W$, such that
  $W'_{|W} = E_0 + E_1$, with $E_0 \geq 0$ and $E_1 \in \Lambda$. By Hironaka's
  theory 
\cite{6}, we may use blowing-ups with smooth centers contained in $W \cap W'$, such that after taking the proper transforms of the elements of our pencil, to get:
  \vskip4pt
  
(a) $  \Supp (E_0)$ has normal crossing;

(b) $\Lambda$ is base-points free (so $\varphi : W \to \P^{n-1}$ is a birational morphism).
\vskip4pt

Further blowing-up of the components of $\Supp (E_0)$ allows to assume $E_0 =
0$ so $D_{|W}$ is linearly equivalent to $E$. Using the previous exact
sequence and the fact that $H^1(X, \O _X) = 0$, it follows that $Bs|D| =
\emptyset$. Finally, $(D^n) = (D_{|W})_W^{n-1} = 1$, so $\varphi$ is a
birational morphism to $\P^n$.
\end{proof}
\begin{example} (cf.\ \cite[Example 1.4]{7}) Let $X \subset \P^{n+d-2}$ be a non-degenerate projective
    variety of dimension $n\geq 2$ and degree $d\geq 3$, which   is not a
    cone.
Then $X$ is rational, unless it is a smooth cubic hypersurface, $n\geq 3$.

If $X$ is singular,  by projecting  from a singular  point we get a variety of
minimal degree, birational to $X$. So $X$ is rational. If $X$ is not linearly
normal, $X$ is isomorphic to a variety of minimal degree. Hence we may assume
$X$ to be smooth and linearly normal. One sees easily that such a linearly
normal, non-degenerate manifold $X\subset \P^{n+d-2}$ has anticanonical
divisor linearly equivalent to $n-1$ times the hyperplane section, i.e.\ they
are exactly the so-called ``classical del Pezzo manifolds''. They were
classified  by Fujita in a series of papers; see \cite{9} for a survey of his
argument.
As Fujita's proof is quite long and difficult, we show how
Theorem~\ref{trei-sase} may be used to prove directly the  rationality of $X$
if $d\geq 4$. Consider the surface $W$ obtained by intersecting $X$ with $n-2$
general hyperplanes. Note that $W$ is a nondegenerate, linearly normal surface 
of degree $d$ in $\P^d$, so it is a del Pezzo surface. As such, $W$ is known
to admit a representation $\varphi:W \to \P^2$ as the blowing-up of $9-d$
points. Let $L\subset W$ be the pull-back via $\varphi$
of a general line in $\P^2$. It is easy to see that $L$ is a cubic rational
curve in the embedding of $W$ into $\P^d$. So, for $d\geq 4$ $L$ is contained
in a hyperplane of $\P^d$. This shows that the conditions of the
Theorem~\ref{trei-sase} are fullfiled for $X$, $|D|$ being the system of
hyperplane sections. We also see that Theorem~\ref{trei-sase} is sharp, as the
previous argument fails exactly for the case of cubics. 
\end{example}
\section{Formal geometry of quasi-lines}\label{S5}
If $Y$ is a closed subscheme of the scheme $X$, the theory of formal functions
of $X$ along $Y$ was developed by Zariski and Grothendieck as an algebraic
substitute for a complex tubular neighbourhood of $Y$ in $X$. We denote by
$X_{/Y}$ the formal completion of $X$ along $Y$, which is the ringed space with
topological space $Y$ and sheaf of rings $\invlim{}_n{\O _{X}/\mathcal{I}^n}$, $\I$ being
the sheaf of ideals defining $Y$ in $X$. All results from formal geometry we
shall need may be found in either R. Hartshorne's classic \cite{4}, or in the
recent comprehensive monograph by L. B\u adescu \cite{1}. Let us recall that
according to Hironaka--Matsumura one can define the ring of formal rational
functions of $X$ along $Y$, denoted by $K(X_{/Y})$. In good cases, it is a
field containing $K(X)$, the field of rational functions of the variety
$X$. We say $Y$ is {\rm G2} in $X$ if $K(X_{/Y})$ is a field and the degree
$[K(X_{/Y}): K(X)]$ is finite. An important result due to Hartshorne implies
that if $Y$ and $X$ are projective manifolds and $N_{Y/X}$ is ample, $Y$ is
{\rm G2} in $X$. We say $Y$ is {\rm G3} in $X$ if $Y$ is {\rm G2} and the inclusion $K(X) \subset K(X_{/Y})$ is an equality. The following definition (cf.\ \cite{7}) is similar to the one in Definition~\ref{trei-patru}; here for a given quasi-line we consider its \' etale neighbourhoods instead of its Zariski neighbourhoods. 
\begin{defin}\label{patru-unu}
 Let $Y \subset X$ be a quasi-line. The number $\widetilde e(X,Y)$ is the
 minimum of $e(X',Y')$, where $X'$ is a projective manifold, $Y'
 \subset X'$ is a quasi-line and $f:(X',Y') \to (X,Y)$ is a morphism of
 models, 
 \' etale along~$Y'$.
\end{defin}

The following theorem, based on results due to Hartshorne and Gieseker, is essential for the sequel.
\begin{theorem}\label{patru-doi} {\rm  (cf.\ \cite[Theorem 5.4]{7})}
If $Y \subset X$ is a quasi-line then $$e(X,Y) = \widetilde{e}(X,Y) \cdot [K(X_{/Y}):K(X)].$$
\end{theorem}
\begin{proof}
 One first observes that (see Lemma~\ref{doi-doi}), if $f:(X',Y') \to
 (X,Y)$ is a morphism of models, with $Y$,$Y'$ quasi-lines, \' etale
 along $Y'$, then we have $e(X,Y) = \deg f\cdot e(X',Y')$. Next we
 choose an $f$ as above, such that $e(X',Y') = \widetilde{e}(X,Y)$ and
 we claim that $Y'$ is {\rm G3} in $X'$. Since $Y'$ is {\rm G2} in
 $X'$, we may apply a very useful construction due to
 Hartshorne--Gieseker (see \cite {1}) to get a
 morphism of models $g:(X'',Y'') \to (X',Y')$ as above, such that $\deg
 g = [K(X'_{/Y'}):K(X')]$ and $Y''$ is {\rm G3} in $X''$. By the previous
 step and the definition of 
$\widetilde{e}(X,Y)$, it follows that $[K(X'_{/Y'}):K(X')] = 1$, i.e. $Y'$ is {\rm G3}~in~$X'$.
 
 The diagram associated to $f$,
\[ \begin{array}{ccc}
K(X)& \longrightarrow& K(X_{/Y})\\[6pt]
\Big\downarrow&&\Big\downarrow\vcenter{\rlap{$\scriptstyle \simeq$}}\\[10pt]
K(X')&\overset{\sim}{\longrightarrow}&K(X'_{/Y'})\end{array}\] 
shows that $\deg f = [K(X_{/Y}):K(X)]$. Note that the right vertical isomorphism comes from the fact that $f$, being \' etale along $Y'$ induces an isomorphism between $X_{/Y}$ and $X'_{/Y'}$.
\end{proof}
\begin{cor} \label{patru-trei}{\rm (cf.\ \cite[Corollary 5.5]{7})}
 Let $(X,Y)$ and $(X',Y')$ 
be models with $Y, Y'$ quasi-lines. If $X_{/Y}\overset{\sim}{\to} X'_{/Y'}$ as formal schemes, then $\widetilde{e}(X,Y)=\widetilde{e}(X',Y')$.
\end{cor}
\begin{cor}\label{patru-patru}{\rm (cf.\ \cite[Corollary 5.6]{7})}
 Let $Y\subset X$ be a quasi-line. Then 
 $e_0(X,Y)\leq \widetilde{e}(X,Y)\leq e(X,Y)$.
\end{cor}
\begin{defin}\label{patru-cinci}
 We say that a quasi-line $Y\subset X$ is {\it regular} if 
$\widetilde{e}(X,Y) = e(X,Y)$.
\end{defin}
\begin{cor}\label{patru-sase}{\rm (cf.\ \cite[Corollary 5.7]{7})}
 Let $Y \subset X$ be a quasi-line. $Y$ is regular if and only if $Y$ is {\rm G3} in $X$. If $e_0(X,Y) = e(X,Y)$, then $Y$ is regular.
\end{cor}

 Note, as a very special case, that a quasi-line $Y \subset X$ with $e(X,Y) =1$ is {\rm G3}. This generalizes the fact, first noticed by Hironaka in the sixties, that a line in $\P^n$ is {\rm G3}.
\begin{cor}\label{patru-sapte}
 Let $(X,Y)$, $(X',Y')$ be models with $Y$, $Y'$ regular quasi-lines. 
If $X_{/Y}\overset{\sim}{\longrightarrow}X'_{/Y'}$ as formal schemes, 
then $(X,Y)\sim (X',Y')$.
\end{cor}

The following proposition generalizes \cite[Corollary 5.12]{7}.
\begin{prop}\label{cinci-*}
 Let $X,X' \subset \P^{d+1}$ be smooth Fano threefolds 
 of degree $d$, $d= 3,4,5$ and let $Y\subset X$, $Y'\subset X'$ be general conics. 
 If $X_{/Y}\overset{\sim}{\longrightarrow}X'_{/Y'}$, then there exists an isomorphism of models $\varphi:(X,Y) \to (X',Y')$.
\end{prop}

For a proof, combine Example~\ref{unu-cincispe} with Corollaries~\ref{patru-sase} and \ref{patru-sapte}.
\begin{defin}\label{patru-zece}
$Y\subset X$ is a {\it line} if $Y$ is a regular almost-line. 
\end{defin}

An ordinary line in $\P^n$ is clearly a line in the sense of Definition~\ref{patru-zece}, whence the terminology.

Using Theorem~\ref{unu-doispe} and the Hartshorne--Gieseker construction one sees that the following holds: 
\vskip4pt

Given a rationally connected  projective manifold $X$, we may find a generically finite morphism $X'\to X$ such that $X'$ contains a line. 
\vskip4pt

The next question looks interesting:
\begin{question}[existence of lines] \label{patru-unspe}
Given a rationally connected manifold $X$ can we find a (smooth, projective) birational model of $X$ containing a line?
\end{question}

Note that almost-lines $Y\subset X$ with $e(X,Y)=1$ are lines by Corollary~\ref{patru-sase}.

The following proposition will allow us to construct other examples of lines.
\begin{prop}\label{patru-doispe}
In the setting of Theorem~{\rm \ref{unu-trei}} assume moreover that $d=1$
and $\widetilde Y$ is {\rm G3} in $Z$. Then $\widetilde Y$ is {\rm G3} in $\widetilde{X}$.
\end{prop}
\begin{proof} Assume that $\widetilde Y$ is not G3 in $\widetilde X$ and apply the 
Hartshorne--Gieseker construction to the model $(\widetilde X, \widetilde Y)$. We find the morphism of models $f: (\widetilde{\widetilde X}, \widetilde {\widetilde Y}) \to (\widetilde X, \widetilde Y)$, \' etale along $\widetilde{\widetilde Y}$ and of degree $>1$. Use Bertini's Theorem to infer that, for a general line $l\subset \P^s$, $\widetilde Z:= (\varphi \circ f)^{-1}(l)= f^{-1}(Z)$ is smooth and irreducible. We find that $f^{-1}(\widetilde Y)\subset \widetilde Z$ is disconnected, having $\widetilde {\widetilde Y}$ as a connected component. A fundamental result due to Hironaka--Matsumura (see \cite{1}) asserts that $f^{-1}(\widetilde Y)$ is G3 in $\widetilde Z$, because $\widetilde Y$ is G3 in $Z$. Now, it is easy to see that a G3 subscheme of a projective manifold is connected (see \cite{1}). This is a contradiction, so $\widetilde Y$ is G3 in $\widetilde X$.   
\end{proof}
\begin{cor}\label{patru-treispe}
 In the setting of Corollary~{\rm \ref{unu-cinci}}, assume moreover that $d=1$. Then $Y$ is {\rm G3} in $X$.
\end{cor}

We may apply the preceding proposition, noting that, $Z$ being a surface, any curve  with positive self-intersection on it is G3.
\begin{example} \label{patru-paispe}
 Let $X$ be the blowing-up of a smooth cubic threefold in $\P^4$ 
with center an ordinary line. $X$ carries a conic-bundle structure $\varphi:X
\to \P^2$; if $Y \subset X$ is a section for $\varphi$ of self-intersection
one, $Y$ is a line. To see this, notice that $Y$ is an almost-line and apply 
Corollary~\ref{patru-treispe}. 
\end{example}

 The following question seems very interesting (especially if it has an affirmative answer), but looks difficult.
\begin{question}\label{patru-cinspe}
 Let $Y\subset X$ be a {\it line}. Is it true that  $e(X,Y)$ is a birational invariant of $X$ (i.e. an invariant of $K(X_{/Y})$)?
\end{question}

The following example shows that the answer is negative if we only assume $Y$
to be a {\it regular quasi-line}.
\begin{example}\label{patru-saispe}
 Let $X \subset \P^5$ be a smooth complete intersection of two quadrics
 and let $Y\subset X$ be a general conic, cf.\ Example~\ref{unu-cincispe}. From Corollary~\ref{patru-sase} we deduce that $Y$ is a regular quasi-line. However, $X$ being rational, a positive answer to Question~\ref{patru-cinspe} would imply $e(X,Y)=1$.
\end{example}

The next result, inspired by Proposition~\ref{doi-unu}, gives a characterization of the
formal completion $\P^n_{/{\rm line}}$.

\begin{theorem}\label{patru-saptispe} 
 Let $(X,Y)$ be a model.
 
\hskip6pt{\rm (i)} There is some $L \in \Pic (X_{/Y})$ such that $\deg L_{|Y}=1$.

\hskip3pt{\rm (ii)} For any such $L$, we have $h^0 (X_{/Y},L)\leq n+1$.

{\rm (iii)} $X_{/Y}$ is isomorphic to $\P^n_{/{\rm line}}$ if and only if there is an $L\in \Pic (X_{/Y})$ such that $\deg L_{|Y}=1$ and $h^0 (X_{/Y},L)= n+1$.
\end{theorem}
\begin {proof}
Denote by $Y(i)$, $i\geq 0$, the i-th infinitesimal neighbourhood of $Y$ in $X$. We have the standard exact sequence 
\[ 0 \to S ^{i+1} ( N^{\scriptscriptstyle\vee}_{Y/X} )
\to \O_{ Y(i+1)} 
\to \O _ { Y(i)} 
\to 0.\]

(i) The above sequence yields the truncated exponential sequence
\[ 0 \to S^{i+1} (N^{\scriptscriptstyle \vee} _{Y/X}) \to \O ^*_{Y(i+1)} \to \O^*_{Y(i)} \to 1.\]
For any curve $Y$, we have $H^2 (Y, S^{i+1} (N^{\scriptscriptstyle\vee}_{Y/X}))=0$, so, by taking cohomology, we get surjections
\[\Pic(Y(i+1))\to \Pic (Y(i))\to 0 \quad \mbox{for } i\geq 0.\]
 Therefore we may lift $\O_{\P^1}(1)$ to $\Pic (X_{/Y})= \invlim_n \Pic(Y(n))$. 
 
 (ii) Let $L_i:= L_{|Y(i)}$ for $i\geq 0$. The first exact sequence above, tensored by $L$, gives
 \[0 \to S^{i+1}(N^{\scriptscriptstyle\vee}_{Y/X})\otimes L_0\to L_{i+1} \to L_i\to 0.\]
 We deduce
 easily
 \[h^0(L_1)\leq n+1, \quad h^0(L_{i+1}) \leq h^0(L_i) \quad \mbox{for } i \geq 1.\]
 As $H^0 (X_{/Y}, L)=\invlim_n H^0(Y(n), L_n)$, the conclusion follows.
 
 (iii) One implication is obvious. To see the other, we remark that the hypothesis $h^0(X_{/Y}, L) =n+1$ and the preceding exact sequences yield that $h^0(Y(i), L_i)= n+1$ for $i\geq 1$, and each $L_i$ is spanned by global sections. Now it is easy to see, using exact sequences as above and induction on $i$, that $L_i$ induces an isomorphism of schemes between $Y(i)$ and the i-th infinitesimal neighbourhood of a line in $\P^n$. These isomorphisms are compatible, so they patch together to give the desired isomorphism $X_{/Y} \overset{\sim
} {\longrightarrow} \P^n_{/{\rm line}}$.
\end{proof}

The next  corollary is the third promised characterization  of the model 
$\Pl$.
\begin{cor}\label{patru-optispe}
 The following conditions are equivalent for a model $(X,Y)$:
 
{\rm (a)} $(X,Y)\sim \Pl$;

{\rm (b)} {\rm (i)} there is some $L\in \Pic(X_{/Y})$ such that 
$\deg L_{|Y}=1$ and $h^{0}(X_{/Y},L) $ $\geq n+1$; 

\hskip15pt{\rm (ii)} $Y$ is regular.
\end{cor}
 
Note that condition (b)(i) implies that $Y$ is a quasi-line, so that (b)(ii) makes sense. To see that (b)(i) and (b)(ii) imply (a), combine Theorem~\ref{patru-saptispe} and Corollary~\ref{patru-sapte}.
\begin{prop}\label{patru-noospe} {\rm (cf.\ \cite[Proposition 4.2]{8}, \cite[Lemma 5.9]{7})}
 Let  $(X,Y)$ be a model with $Y$ a quasi-line. Let $E$ be a vector
 bundle on $X$ such that 
\[E_{|Y}= \O_{\P^1}(a)\oplus \cdots \oplus \O_{\P^1}(a) \oplus \O_{\P^1} (a+1)\] 
for some $a\in \Z$. Let $X'$ be  $\P(E)$ and 
let $\pi:X'\to X$ be the projection. Then:
 
\hskip6pt{\rm (i)} there is a quasi-line $Y' \subset X'$ such that
 $\pi:(X',Y')\to (X,Y)$ is a morphism of models;

\hskip3pt{\rm (ii)} $e(X,Y)=e(X',Y')$.
\end{prop}
\begin{example}\label{patru-douazeci} 
{\rm (i)} The model $(X,Y)$ from Example~\ref{unu-paispe} satisfies
condition (b)(i) from Corollary~\ref{patru-optispe}, but $Y$ is not regular: we have $\widetilde{e}(X,Y)=1$ and  $e(X,Y)=n+1$.

{\rm (ii)} Consider  the model $\Pl$ and apply Proposition~\ref{patru-noospe}
to $E=T_{\P^n}$ to find the new model $(X'=\P(T_{\P^n}),Y')$. $(X',Y')$
satisfies (b)(ii) of Corollary~\ref{patru-optispe}, but does not satisfy
(b)(i). Indeed $Y'$ is regular since $e(X',Y')=e\Pl=1$ by
Proposition~\ref{patru-noospe}(ii). 
If (b)(i) would be fulfilled, we would have $(X',Y')\sim (\P^{2n-1}, \mbox{\rm
  line})$. 
But $(X',Y')$ is easily seen to be minimal (see \cite[Lemma 4.4]{8}), so $(X',Y')$ would be isomorphic to $(\P^{2n-1}, \mbox{line})$, which is clearly absurd. Note that $(X',Y')$ also provides an example verifying the first two conditions of Theorem~\ref{doi-cinci}(b), but not the third.
\end{example}
\section*{Acknowledgements} We thank the organizers of the Siena
Conference for   inviting us. Their dedication resulted in a very
enjoyable mathematical and social encounter, in one of the most
beautiful surroundings one can dream of.

We benefitted from two one-month visits, first at the University of
Genova in September 2003 and the other at Universit\' e  Lille 1 in April
2004, while working on this paper. We thank Lucian B\u adescu, Mauro 
Beltrametti and Jean d'Almeida for their kind invitations and for
making our stay very enjoyable.

We acknowledge partial financial support from Contract CNCSIS\linebreak
no. 33079/2004.

\end{document}